\documentclass[12pt]{article}

\usepackage{latexsym}
\usepackage{amsfonts}
\usepackage{amssymb}
\usepackage{amsmath}
\usepackage{mathrsfs}
\usepackage{bbm}
\usepackage{dsfont}
\usepackage{mathtools}
\DeclarePairedDelimiter\norm{\lVert}{\rVert}
\usepackage{enumitem}
\usepackage{tikz}
\usepackage{color}
\definecolor{indigo}{RGB}{0,0,120}
\usepackage[colorlinks=true,linkcolor=indigo,citecolor=blue,urlcolor=indigo]{hyperref}
\usepackage{lettrine}

\newcommand{\be}{\begin{equation}}
\newcommand{\ee}{\end{equation}}
\newcommand{\bea}{\begin{eqnarray}}
\newcommand{\eea}{\end{eqnarray}}
\newcommand{\bean}{\begin{eqnarray*}}
\newcommand{\eean}{\end{eqnarray*}}
\newcommand{\brray}{\begin{array}}
\newcommand{\erray}{\end{array}}

\newtheorem{dfn}{Definition}[section]
\newtheorem{thm}[dfn]{Theorem}
\newtheorem{lmma}[dfn]{Lemma}
\newtheorem{ppsn}[dfn]{Proposition}
\newtheorem{crlre}[dfn]{Corollary}
\newtheorem{xmpl}[dfn]{Example}
\newtheorem{rmrk}[dfn]{Remark}

\newcommand{\bdfn}{\begin{dfn}\rm}
\newcommand{\bthm}{\begin{thm}}
\newcommand{\blmma}{\begin{lmma}}
\newcommand{\bppsn}{\begin{ppsn}}
\newcommand{\bcrlre}{\begin{crlre}}
\newcommand{\bxmpl}{\begin{xmpl}}
\newcommand{\brmrk}{\begin{rmrk}\rm}

\newcommand{\edfn}{\end{dfn}}
\newcommand{\ethm}{\end{thm}}
\newcommand{\elmma}{\end{lmma}}
\newcommand{\eppsn}{\end{ppsn}}
\newcommand{\ecrlre}{\end{crlre}}
\newcommand{\exmpl}{\end{xmpl}}
\newcommand{\ermrk}{\end{rmrk}}

\title{Decomposability of multiparameter CAR flows}
\author{Anbu Arjunan}

\begin{document}
\maketitle
\begin{abstract}
 Let $P$ be a closed convex cone in $\mathbb{R}^d$ which is assumed to be spanning $\mathbb{R}^d$ and contains no line. In this article, we consider a family of CAR flows over $P$ and study the decomposability of the associated product systems. We establish a necessary and sufficient condition for CAR flow to be decomposable. As a consequence we show that there are uncountable many CAR flows which are cocycle conjugate to the corresponding CCR flows.  
 \end{abstract}
\noindent {\bf AMS Classification No. :} {Primary 46L55; Secondary 46L99.}  
\section{Introduction}
We study the decomposability of CAR flows over closed convex cones.
In the one-parameter case, there exists only countable many CCR flows and they exhaust all type I $E_0$-semigroups; see \cite{Arv}. It is known that CCR flows and CAR flows are cocycle conjugate \cite{PC}. But in the multiparameter context there are uncountable many CCR flows over a closed convex cone $P$ \cite{ASS,AS19,sundar2019arvesons}. It was shown in \cite{R19} that for certain $P$-modules, CCR flow is not cocycle conjugate to the CAR flow. 
Let $P$ be a closed convex cone in $\mathbb{R}^d$ and let $V$ be an isometric representation of $P$. 
Denote by $\alpha$ the CAR flow associated to the isometric representation $V$ and by $E$ the product system associated to $\alpha$.
A natural question that arises is the following. Under what condition on the isometric representation $V$, the corresponding product system $E$ is decomposable $?$. In this paper we answer this for the isometric representations given by $P$-modules.  More precisely, we show that the CAR flow associated to a P-module is decomposable if and only the module is a half-space.   

The organization of the paper is as follows.
In section 1. we establish a bijective correspondence between the set of all additive decomposable sections and the set of all coherent sections for the product system over $\mathbb{R}_+$ which is assumed to have a coherent section.
In section 2. For a fixed $P$-module and a ray in a closed convex cone $P$, by using the correspondence obtained in section 1, we compute the decomposable vectors of the CCR flow along the ray.
In section 3. We show that there are uncountable many CCR flows which are cocycle conjugate to CAR flows over a closed convex cone.
\section{Non-commutative stochastic calculus}
The sole purpose of this section is to record for future reference that the bijection between addits and units established in \cite{OR} works equally well to provide a bijection between the set of all coherent sections and the space of all centred additive decomposable sections and we do not claim much originality.  As an immediate application, we obtain another proof for the $e$-Logarithm of a coherent section $e=(e_t)_{t\in \mathbb{R}_+}$ is positive definite.
We leave it to the reader for proof of many results and we only provide the proof for the necessary places. We also adapt most of the notations from \cite{OR}. 

Let $E=\{(t,\xi): t\geq 0, \; \xi\in E_t\}$ be a product system and in short we write $E=\bigcup_{t\geq 0}E_t$.
A family $\{x_s:0\leq s<t\}$ of vectors is said to be \emph{left coherent} if each $x_s\in E_s$ and for $0\leq r<s<t$, there exists a vector $x(r,s)\in E_{s-r}$ such that $x_s=x_rx(r,s).$
A \emph{left coherent section} of $E$ is a left coherent family $\{x_s:0<s<t\}$ with $t=\infty$ (similarly we also have a notion of right coherent section). 
Now onwards we simply call the left coherent family a coherent family.
Then the family $\{x(r,s):0\leq r< s<t\}$ satisfies 
$x(q,r)x(r,s)=x(q,s)\text{ for each } 0\leq q<r<s.$ 
Now onwards, we assume that all the product system in our discussion will have a coherent section.
Fix a coherent section $\Omega=(\Omega_t)_{t\geq 0}$ along with the product system $E$ such that $\|\Omega_t\|=1$ for $t\geq 0$.
For a fixed $r\in \mathbb{R}_+$, let $\widetilde{\Omega}_{t}=\Omega(r,r+t)$ for $t\geq 0.$ Note that $\widetilde{\Omega}=(\widetilde{\Omega}_t)_{t\in\mathbb{R}_+}$ defines a coherent section.
A coherent section $(x_t)_{t\in \mathbb{R}_+}$ is called a centred coherent section with respect to $\Omega$ (or simply a centred coherent section) if $\langle x_t|\Omega_t\rangle =1$ for each $t\geq 0$.
\begin{dfn}\label{additive_decomposible}
Let $\Omega=\{\Omega_s:0\leq s<t\}$ be a coherent family.
 A family $\{b_s: 0\leq s<t\}$ of vectors is said to be left additive decomposable with respect to $\Omega$ if for any $0\leq r<s<t,$ there exists a vector $b(r,s)\in E_{s-r}$ such that
 \[b_r\otimes \Omega(r,s)+ \Omega_r\otimes b(r,s)=b_s.\]
 We simply say $(b_t)_{t\in\mathbb{R}_+}$ as additive decomposable section if $\Omega$ is clear from the context.
 If in addition the family $\{b_s: 0\leq s<t\}$ satisfies $\langle b_s| \Omega_s\rangle =0$ for each $0\leq s<t.$ We say that the family $\{b_s:0\leq s<t\}$ a centred additive decomposable family.
\end{dfn}
Every additive decomposable family $(b_t)_{t\in\mathbb{R}_+}$ has the following decomposition: For $t\geq 0$, $b_t=c_t+\langle b_t|\Omega_t\rangle \Omega_t$ where $(c_t)_{t\in\mathbb{R}_+}$ is a unique centred additive decomposable family.
We observe that the family $\{b(r,s):0\leq r< s<t\}$ satisfies $b(q,s)=b(q,r)\otimes \Omega(r,s)+\Omega(q,r)\otimes b(r,s)$ for each $0\leq q<r<s.$ 
If an additive decomposable family $\{b_s:0\leq s<t\}$ is centred, we also have $\langle b(r,s)| \Omega(r,s)\rangle=0 $ for every $0\leq r<s<t.$ 
An \emph{additive decomposable section} is an additive decomposable family $\{b_s:0\leq s<t\}$ with $t=\infty$. 
\begin{lmma} \label{limit_zero}
Let $\Omega=\{\Omega_s:0\leq s<t\}$ be a left coherent family with $\|\Omega_s\|=1$ for each $s$ and 
let $b=\{b_s:0\leq s<t \}$ and $c=\{c_s:0\leq s<t \}$ be centred additive decomposable families with respect to $\Omega$. Then $\displaystyle\lim_{h\to 0^+}\langle b_h|c_h\rangle=0.$ 
\end{lmma}
\begin{prf}
 For $0<r<s$, define a map $L_{r}:E_{s-r}\to E_{s}$ by $L_r(\xi)=\Omega_r\xi$ for $\xi \in E_{s-r}$.
 Then $L_r$ is an isometry. Denote the range projection $L_rL_r^*$ of $L_r$ by $P_r$  i.e. we can view $P_r=\langle .|\Omega_r\rangle\Omega_r \otimes 1_{E_{s-r}}:E_r\otimes E_{s-r}\to E_r\otimes E_{s-r}$.
 Then $P_r$ is strongly converges to $1$ as $r$ decreases to zero (see \cite[Theorem 6.1.1]{Arv} for the proof of this fact).
Let $0<h<s<t$. We observe that $\|b_s\|^2=\|b_h\|^2+\|b(h,s)\|^2$ and we have
\begin{align*}
 \|b_s\|^2 &=\displaystyle \lim_{h\to 0^+}\| P_h(b_s)\|^2\\
           &=\displaystyle \lim_{h\to 0^+}\|P_h(b_h\otimes \Omega(h,s)+\Omega_h \otimes b(h,s) )\|\\
           &=\displaystyle \lim_{h\to 0^+}\|b(h,s) )\|\\
           &=\displaystyle \lim_{h\to 0^+}\|b_s\|^2-\|b_h\|^2.
 \end{align*}
From the above, we conclude that $\displaystyle\lim_{h\to 0^+}\| b_h\|^2=0$ and hence $\displaystyle\lim_{h\to 0^+}\langle b_h|c_h\rangle=0$.  \hfill $\Box$
 \end{prf}
\begin{ppsn} \label{continuity_of_additive_vector}
 Let $\{b_s:0<s<t\}$ and $\{c_s:0<s<t\}$ be centred additive decomposable families. Then the map
 $(0,t)\ni s\mapsto \langle b_s | c_s\rangle \in \mathbb{C}$ is continuous.
\end{ppsn}
\begin{prf}
 Let $s_0\in (0,t)$ be given. We show that the map $ (0,t)\ni s\mapsto \langle b_s |c_s\rangle \in \mathbb{C}$ is both left and right continuous.
 For any $0<h<t-s_0$, we observe that
 \begin{align*}
  \langle  b_{s_0+h}|c_{s_0+h} \rangle&=\langle  b_{s_0}|c_{s_0} \rangle+\langle  b(s_0,s_0+h)|c(s_0,s_0+h) \rangle\\
        &=\langle  b_{s_0}|c_{s_0} \rangle+\langle \widetilde{ b}_h|\widetilde{c}_h \rangle \text{ where }\widetilde{b}_{h}=b(s_0,s_0+h)\text{ and }\widetilde{c}_{h}=c(s_0,s_0+h) .
 \end{align*}
Then the sets $\{\widetilde{b}_{r}:0<r<t-s_0\}$ and $\{\widetilde{c}_{r}:0<r<t-s_0\}$ form centred left additive decomposable families.
Indeed, for $0<q<r<t-s_0$, we have
\[\widetilde{b}_r=\widetilde{b}_{q}\otimes \widetilde{\Omega}(q,r)+\widetilde{\Omega}_{q}\otimes \widetilde{b}(q,r)\] where $\widetilde{\Omega}(q,r)=\Omega(s_0+q,s_0+r)$, $\widetilde{\Omega}_q=\Omega(s_0,s_0+q)$, and $\widetilde{b}(q,r)=b(s_0+q,s_0+r)$.
Similarly $\{\widetilde{c}_{r}:0<r<t-s_0\}$ forms a centred left additive decomposable family.
By Lemma \ref{limit_zero}, $\displaystyle\lim_{h\to 0^+}\langle \widetilde{ b}_h |\widetilde{ c}_h\rangle =0$. Hence the equation $\langle  b_{s_0+h}|c_{s_0+h} \rangle =\langle  b_{s_0}|c_{s_0} \rangle+\langle \widetilde{ b}_h|\widetilde{c}_h\rangle$ implies that $\langle  b_{s_0+h}|c_{s_0+h} \rangle \to\langle  b_{s_0}|c_{s_0} \rangle$ as $h$ goes to zero. This means that the map $(0,t)\ni s\mapsto \langle b_s|c_s\rangle \in \mathbb{C}$ is right continuous.
\\
For $0<h<s_0$, let $b':=\{b'_s=b(s_0-s,s_0): 0<s<s_0\}$ and $\Omega':=\{\Omega'_s=\Omega(s_0-s,s_0): 0<s<s_0\}$.
Note that $\Omega'$ is a right coherent family. We leave it to the reader to verify that $b'$ is a centred right additive decomposable family with respect to $\Omega'$ for the product system $E$.
In the opposite product system $E^{\text{op}}$, $\Omega'$ is a left coherent family and $b'$ is a centred left additive decomposable family with respect to $\Omega'$.
Similar to the above argument,
we see that $\langle  b_{s_0-h}|c_{s_0-h} \rangle\to \langle b_{s_0}|c_{s_0} \rangle$ as $h$ goes to zero. Hence the $(0,t)\ni s\mapsto \langle b_s|c_s\rangle \in \mathbb{C}$ is left continuous.
This proves the proposition. \hfill $\Box$
\end{prf}
\begin{dfn}
\begin{enumerate}
\item By an \emph{adapted process} of $E,$ we mean a measurable map $\mathbb{R}_{+}\ni t\mapsto x_t\in  E$ satisfying $x_t\in E_t$ for each $t\geq 0.$ We say an adapted process $(x_t)_{t\in\mathbb{R}_+}$ is simple if there exists a partition $0\leq t_0<t_1<\cdots <t_n<\cdots $ of $\mathbb{R}_{+}$ such that 
\begin{equation} \label{simple-adapted-equation}
 x_t=\sum_{i=0}^{\infty} x_i \otimes \Omega(t_i,t) \;\chi_{[t_i,t_{i+1})}(t), \text{ for each }t\geq 0\text{ and } x_i\in E_{t_i} \text{ for }i\geq 0 .
\end{equation}
\item Let $(x_t)_{t\in\mathbb{R}_+}$ be a simple adapted process of the form (\ref{simple-adapted-equation}) and let $(b_t)_{t\in\mathbb{R}_+}$ be a centred additive decomposable section of $E$. Then the integral of a simple adapted process $(x_t)_{t\in\mathbb{R}_+}$ with respect to the centred additive decomposable section $(b_t)_{t\in\mathbb{R}_+}$ over the interval [a,b] is denoted by $\int_{a}^{b} x_t db_{t}$ and is defined as follows:
 \begin{equation}\label{stochastic_integral}
  \int_{a}^{b} x_t db_{t}=\sum_{i=m}^{n-1} x_i\otimes b(t_i, t_{i+1})\otimes \Omega(t_{i+1},b).
 \end{equation}
Here we have refined the partion such that $a=t_m$ and $b=t_n$ for some $m,n\in \mathbb{N}$.
\item Let $(x_t)_{t\in\mathbb{R}_+}$ be an adapted process and let $(b_t)_{t\in\mathbb{R}_+}$ be a centred additive decomposable section of $E$. We say that an adapted process $(x_t)_{t\in\mathbb{R}_+}$ is It\^o integrable with respect to $(b_t)_{t\in\mathbb{R}_+}$ if there exists a sequence of simple adapted process $\{x^{(n)}\}_{n=1}^{\infty}$ with $x^{(n)}=(x^{(n)}_t)_{t\in\mathbb{R}_+}$ for each $n\in \mathbb{N}$ such that the sequence $\int_{a}^{b} x^{(n)}_t db_{t}$ is Cauchy for each $a,b\in \mathbb{R}_+$ with $a<b$. In that case we define $\int_{a}^{b} x_t db_{t}$ as follows:
\begin{equation} \label{stochastic_limit}
  \int_{a}^{b} x_t db_{t}:=\lim_{n\to \infty}\int_{a}^{b} x^{(n)}_t db_{t}.
 \end{equation}
 We will see that the above integral is well defined.
\end{enumerate}
\end{dfn}
The definition for It\^o integrable adapted process is slightly different from the one considered in \cite[see discussion after Proposition 5.4]{OR}.
\begin{lmma}\label{basic_properties}
 Let  $(x_t)_{t\in\mathbb{R}_+}$ and $(y_t)_{t\in\mathbb{R}_+}$ be It\^o integrable adapted processes with respect to the centred additive decomposable section $(b_t)_{t\in\mathbb{R}_+}$. Then we have the following properties.
 \begin{itemize}
  \item[(1)]  For $0\leq s\leq t$, we have $\left\langle\displaystyle \int_{s}^{t} x_r db_r| \Omega_t\right\rangle=0$.
  \item[(2)] For $0\leq s_0\leq s$ and $t\geq 0$,  $\displaystyle\int_{s_0}^{s+t}x_r\;db_r=\int_{s_0}^{s}x_r\;db_r \otimes \Omega(s,s+t)+ \int_{s}^{s+t}x_r\;db_r$.
  \item[(3)] For $s,t\geq 0$, $\displaystyle\int_{s}^{s+t} x_s\otimes y_{r-s} db_r=x_s\otimes \int_{0}^{t} y_r d\tilde{b}_r$.\\ Here $\tilde{b}=(\tilde{b_r})_{r\geq 0}:=(b(s,s+r))_{r\geq 0}$ is a centred additive decomposable section with respect to the coherent section $\widetilde{\Omega}=(\Omega(s,s+r))_{r\geq 0}$.
  \item[(4)] For $0\leq s_0\leq s$ and $0\leq t_0\leq  t$, we have $\left\langle \displaystyle\int_{s+t_0}^{s+t} x_r db_r\Bigg| \int_{s_0}^{s} x_r db_r \otimes \Omega(s,s+t)\right\rangle=0.$
  \item[(5)] For $s,t\geq 0$,  $ \displaystyle\int_{s}^{s+t} \Omega_r db_r=b_{s+t}-b_s\otimes \Omega(s,s+t).$
 \end{itemize}
\end{lmma}
\begin{prf}
 By definition it is enough to prove the above results for a simple adapted process and we leave it to reader's verification.\hfill $\Box$
\end{prf}
\\
Let $(b_t)_{t\in \mathbb{R}_+}$ be a centred additive decomposable section. Define a map $F:\mathbb{R}\to \mathbb{R}$ by 
\[F(s):=\begin{cases} 
         \|b_s\|^2 & \mbox{ if } s>0,\cr
         0 & \mbox{ if } s\leq 0.
        \end{cases}
        \]
        Since $F$ is a non-decreasing right continuous function on $\mathbb{R}$, there exists a unique Borel measure $\mu$ on $\mathbb{R}$ such that $\mu((-\infty,s])=F(s)$ for $s\in \mathbb{R}$. 
        For simple adapted processes  $(x_t)_{t\in \mathbb{R}_+}$ and $(y_t)_{t\in \mathbb{R}_+}$, we have
 \begin{equation}
  \left\langle \int_{s}^{t} x_r db_r\Bigg| \int_{s}^{t} y_r db_r\right\rangle=\int_{s}^{t} \langle x_r|y_r \rangle d\mu(r).
 \end{equation}
 We can see using the above equality that (\ref{stochastic_limit}) is well defined. 
        With this results, we have the following lemma which can also be thought of as a version of It\^o identity.
\begin{lmma} \label{ito_identity}
 Let  $(x_t)_{t\in \mathbb{R}_+}$ and $(y_t)_{t\in \mathbb{R}_+}$ be continuous adapted processes and let $(b_t)_{t\in \mathbb{R}_+}$ be a centred additive decomposable section. Then we have
 \begin{equation}
  \left\langle \int_{s}^{t} x_r db_r\Bigg| \int_{s}^{t} y_r db_r\right\rangle=\int_{s}^{t} \langle x_r|y_r \rangle d\mu(r).
 \end{equation}
\end{lmma}
The following proposition provides a condition for a continuous adapted process to be It\^o integrable and a proof follows from Lemma \ref{ito_identity}.
\begin{ppsn} \label{condition_ito_integral}
Let $(b_s)_{s\in \mathbb{R}_+}$ be a centred additive decomposable section and let $x=(x_t)_{t\in \mathbb{R}_+}$ be a continuous adapted process such that 
\begin{equation} \label{ito_equation}
 \langle x_{r+s}|x_r\otimes \Omega(r,r+s)\rangle =\|x_r\|^2 \text{ for each }r,s\geq 0.
\end{equation}
Then $ x=(x_t)_{t\in \mathbb{R}_+}$ is It\^o integrable with respect to $(b_s)_{s\in \mathbb{R}_+}$. 
In fact the  sequence $x^{(n)}$ of simple adapted process given by 
 $x^{(n)}_r=\sum_{i=0}^{n-1} x_{r^{(n)}_i} \otimes \Omega(r^{(n)}_i,r)\chi_{[r^{(n)}_i,r^{(n)}_{i+1})}(r)$
 for each $r\geq 0$ and $r^{(n)}_i=s+(t-s)\frac{i}{n-1}$ with $0\leq i\leq n-1,$ converges to $ x=(x_t)_{t\geq 0}$. 
\end{ppsn}
\begin{lmma} \label{integration_by_parts_1}
 Let $(b_t)_{t\in \mathbb{R}_+}$ be a centred additive decomposable section and $\mu$ be its associated measure given after Lemma \ref{basic_properties}. Then we have, \[\displaystyle \int_{0}^t \|b_r\|^{2n} d\mu(r)=\frac{ \|b_t\|^{2(n+1)}}{n+1} \;\;\text{    for every }n\geq 0.\]
\end{lmma}
\begin{prf}
Recall that \cite[Theorem 6.5.10]{mstoll} if $f$ and $\alpha$ are continuous monotone nondecreasing functions on $[a,b]$. Then $f\in R(\alpha)$ and $ \alpha\in R(f)$ (here $R(\alpha)$ denotes the space of Riemann-Stieltjes functions with respect to $\alpha$ ). Moreover we have, 
\begin{equation}\label{integration_by_parts}
\int_{a}^{b}f(s)\;d\alpha(s)+\int_a^b\alpha(s)\;df(s)=\alpha(b)f(b)-\alpha(a)f(a).
\end{equation}
Taking $\alpha(s)=f(s)=\|b_s\|^2$. Then for $t>0,$ we have $\int_0^t\|b_s\|^2\; d\mu(s)=\frac{\|b_t\|^4}{2}.$ 
Let us take $\alpha(s)=\|b_s\|^{2(k+1)}$ for $k>1$ and $f(s)=\|b_s\|^2$.
Assume that $ \int_0^t\|b_s\|^{2k}\; d\mu(s)=\frac{\|b_t\|^{2(k+1)}}{k+1} $ is true for $k$. This means that $d\alpha(s)=(k+1)\|b_s\|^{2k}\;d\mu(s)$. 
Then we have
\begin{align*}
 \|b_t\|^{2(k+2)}&= \int_0^t \|b_s\|^2\; d\alpha(s) + \int_0^t \|b_s\|^{2(k+1)} \;d\mu(s) \text{ (by integration by parts) }\\
                 &= \int_0^t \|b_s\|^{2} \; (k+1)\|b_s\|^{2k}\;\;d\mu(s)+\int_0^t \|b_s\|^{2(k+1)}\; d\mu(s)\\
                 &=\int_0^t  (k+2)\|b_s\|^{2(k+1)}\;\;d\mu(s).
\end{align*}
Hence we have shown the lemma by induction. \hfill $\Box$
\end{prf}
\\
Denote the space of all centred additive decomposable sections by $\mathcal{A}$ and the set of all coherent sections by $\mathcal{C}$. The following proposition is very similar to \cite[Proposition 5.9]{OR} we merely provide the sketch of a proof.
\begin{ppsn} \label{construction_exponential}
 Let $(b_t)_{t\in \mathbb{R}_+}$ be a centred additive decomposable section. Then there exists a unique solution to the quantum stochastic integral equation
 \begin{equation} \label{QSIE}
  u_t=\Omega_t+\int_0^tu_s \;db_s \; \text{ for each }t\geq 0, 
 \end{equation}
and the solution is a centred coherent section. Moreover, the map $\mathcal{A}\ni(b_t)_{t\in \mathbb{R}_+}\mapsto (u_t)_{t\in \mathbb{R}_+} \in \mathcal{C}$ is injective.
\end{ppsn}
\begin{prf}
Let $x_t^{(0)}=\Omega_t, \; x_t^{(1)}=b_t, $ and $x_t^{(n)}=\displaystyle\int_{0}^t x_r^{(n-1)}\;db_r$ for $n\geq 1$.  
Define $u_t=\displaystyle\sum_{n=0}^{\infty} x_t^{(n)}$ for $t\geq 0$. 
Then we can see the following.
\begin{enumerate}
 \item $\langle x_t^{(n)}|x_t^{(m)}\rangle =\delta_{m,n}\frac{\|b_t\|^{2n}}{n!}$ and $\|u_t\|^2=e^{\|b_t\|^2}$.
 \item $(u_t)_{t\in \mathbb{R}_+}$ is It\^o integrable and satisfies the quantum stochastic integral equation (\ref{QSIE}). Also by Lemma \ref{integration_by_parts_1}, the solution $(u_t)_{t\in \mathbb{R}_+}$ is unique.
\end{enumerate}
Define $w=(w_t)_{t\in \mathbb{R}_+}$ as 
\[w_r:=\begin{cases} 
         u_r & \mbox{ if } r\in (0,s),\cr
         u_s\otimes \widetilde{ u}_{r-s} & \mbox{ if } r\geq s,
        \end{cases}
\]
where $\widetilde{ u}_{r}=u(s,s+r)$.
Then we see that $w=(w_t)_{t\in \mathbb{R}_+}$ satisfies (\ref{QSIE}). 
By the uniqueness of the solution, $u_{s+t}=w_{s+t}=u_s\otimes \widetilde{ u}_t=u_s\otimes u(s,s+t).$ It is also centred.  
Hence $u=(u_t)_{t\in \mathbb{R}_+}$ is a centred coherent section.
One can check that the map $\mathcal{A}\ni(b_t)_{t\in \mathbb{R}_+}\mapsto (u_t)_{t\in \mathbb{R}_+} \in \mathcal{C}$ is injective.
\hfill $\Box$
\end{prf}
\\
We will use the notation $\text{Exp}(b)$ to denote the solution of (\ref{QSIE}).
Let $(x_t)_{t\in \mathbb{R}_+}$ be any adapted process. From now onwards we fix the following convention. 
\begin{align*}
 x_t'&=x_t\otimes \Omega(t,T)\text{ for any } 0\leq t\leq T, \text{ and}\\
 x_{s,s+t}'&=\Omega_{s}\otimes x_t\otimes \Omega(s+t,T)\text{ for } s,t\geq 0 \text{ with }s+t\leq T.
 \end{align*}
We also provide the sketch for the following proposition.
\begin{ppsn}
 Let $u=(u_t)_{t\in \mathbb{R}_+}$ be a centred coherent section and for $t\geq 0$, $n \in \mathbb{N}$, define 
 $ y_t^{(n)}:=\displaystyle\sum_{i=1}^{2^n} y_t^{i,n}$, where for $1\leq i\leq 2^n$,
 \begin{align*}
  y^{i,n}_t=\Omega_{\frac{t}{2^n}} \otimes \Omega\bigg(\frac{t}{2^n},\frac{2t}{2^n}\bigg)&\otimes\cdots \otimes \Omega\bigg(\frac{(i-2)t}{2^n},\frac{(i-1)t}{2^n} \bigg)\otimes \\
  &\bigg( u\bigg(\frac{(i-1)t}{2^n},\frac{it}{2^n} \bigg)-\Omega\bigg(\frac{(i-1)t}{2^n},\frac{it}{2^n} \bigg) \bigg) \otimes \cdots \otimes \Omega\bigg(t-\frac{t}{2^n},t\bigg).\\
 \end{align*}
Then $\displaystyle \lim_{n\to \infty} y_t^{(n)}$ exists and denote its limit by $\text{Log}(u)_t.$ 
Show that $b=(b_t)_{t\in \mathbb{R}_+}:=(\text{Log}(u)_t)_{t\in \mathbb{R}_+}$ is a centred additive decomposable section. 
\end{ppsn}
\begin{prf}
 For $s>0$, $\|u_s-\Omega_s\|^2=e^{\varphi(s)}-1$ where $\varphi(s)=\|b_s\|^2.$
 Note that $\langle y^{(n)}_t| y^{(m)}_t\rangle =\sum_{k=1}^{2^m}\big( e^{\varphi(\frac{kt}{2^m})-\varphi(\frac{(k-1)t}{2^m})}-1\big) \text{  for }n\leq m$
 and $\|y_t^{(n)}\|^2\to \varphi(t)$ as $n\to \infty$. Now for $n\leq m$, we have
 \begin{align*}
  \|y_t^{(n)}-y_t^{(m)}\|^2&=\|y_t^{(n)}\|^2+\|y_t^{(m)}\|^2-2\text{ Re }\langle y_t^{(n)}|y_t^{(m)}\rangle\\
                           &\to \varphi(t)+\varphi(t)-2\varphi(t)=0 \text{ as } n,m\to \infty.
 \end{align*}
Hence $y_t^{(n)}$ is a Cauchy sequence in $E_t$ and it is convergent.
\\
For $s,t\geq 0$, let $y^{(n)}_{s,s+t}=\sum_{i=1}^{2^n} y^{i,n}_{s,s+t},$ where  
\begin{align*}
  y^{i,n}_{s,s+t}=&\Omega(s,s+\frac{t}{2^n}) \otimes \Omega(s+\frac{t}{2^n},s+\frac{2t}{2^n})\otimes\cdots \otimes \Omega(s+\frac{(i-2)t}{2^n}, s+\frac{(i-1)t}{2^n} )\otimes \\
  &\bigg( u(s+\frac{(i-1)t}{2^n},s+\frac{it}{2^n} )-\Omega(s+\frac{(i-1)t}{2^n},s+\frac{it}{2^n} ) \bigg) \otimes \cdots \otimes \Omega(s+t-\frac{t}{2^n},s+t).\\
 \end{align*}
 Let $\widetilde{\Omega}_t=\Omega(s,s+t)$ and $\widetilde{u}_t=u(s,s+t)$ for each $t\geq 0$.
 Observe that $y^{(n)}_{s,s+t}=\widetilde{y}_{t}^{(n)}$ where $\widetilde{y}_{t}^{(n)}$ is defined using the left coherent section $(\widetilde{u}_t)_{t\in \mathbb{R}_+}$ and the fixed coherent section $(\widetilde{\Omega}_t)_{t\in \mathbb{R}_+}$ similar to the above construction.
 This implies that $\lim_{n\to \infty} y^{(n)}_{s,s+t}$ exists and denote its limit by $b(s,s+t)$.
 First, we claim that for $m\in \mathbb{N}$ and $t\geq 0$, $b'_{mt}=b'_{t}+b'_{t,2t}+\cdots + b'_{(m-1)t,mt}.$
 For it is enough to show that $\|y^{(n)'}_{mt}-\sum_{k=0}^{m-1}y^{(n)'}_{kt,(k+1)t}\|^2\to 0$ as $n\to \infty.$
 Now we observe that
\begin{align*}
  \norm[\bigg]{ y^{(n)'}_{mt}-\sum_{k=0}^{m-1} y^{(n)'}_{kt,(k+1)t}}^2&= \displaystyle\sum_{j=1}^{2^n}\big( e^{\varphi(\frac{jmt}{2^n})-\varphi(\frac{(j-1)mt}{2^n})}-1\big)\\
              &\;\;\quad\quad-\sum_{k=0}^{m-1}\sum_{i=1}^{2^n} \big(e^{\varphi(kt+\frac{it}{2^n})-\varphi(kt+\frac{(i-1)t}{2^n})}-1\big)\\
              &\to \varphi(mt)-\sum_{k=0}^{m-1}\big(\varphi((k+1)t)-\varphi(kt) \big)=0,\text{ as }n\to \infty.
\end{align*}
Hence $b'_{mt}=b'_{t}+b'_{t,2t}+\cdots + b'_{(m-1)t,mt}$, for $n\in \mathbb{N}.$ From which we deduce that $b'_{(m+n)t}=b'_{nt}+b'_{nt,(m+n)t}$ for $0\leq m+n\leq T$.
For any $q\in \mathbb{Q}_+$ with $s+qs\in [0,T]$, we observe that $b'_{s+qs}=b'_{s}+b'_{s,s+qs}$ and $b'_{s,s+t}$ is a limit of $\sum_{i=1}^{2^n} \big(u'_{s+\frac{(i-1)t}{2^n}, s+\frac{it}{2^n}} -\Omega_T \big).$ 
We see that for any $s,t\geq 0$, $b'_{s+t}=b'_{s}+b'_{s,s+t}$.
This implies that $b_{s+t}=b_{s}\otimes\Omega(s,t) +\Omega_s\otimes b(s,s+t)$.
By definition, it is clear that $(b_t)_{t\in \mathbb{R}_+}$ is centred. Hence  $(b_t)_{t\in \mathbb{R}_+}$ is a centred additive decomposable section.  \hfill $\Box$
\end{prf} 
\begin{lmma} \label{lemma_for_theorem}
  Let $(u_t)_{t\in \mathbb{R}_+}$ be a centred coherent section and let $(x_t)_{t\in \mathbb{R}_+}$ be an adapted process satisfying $\langle x_{r+s}|x_r\otimes \Omega(r,r+s)\rangle =\|x_r\|^2 \text{ for each }r,s\geq 0.$ Then $(x_t)_{t\in \mathbb{R}_+}$ is It\^o integrable adapted process with respect to $(\text{Log}_{\Omega}(u)_t)_{t\in\mathbb{R}_+} $. Moreover we have,
 \[\bigg\langle u_t\bigg|\int_{0}^{t} x_r \;d\text{Log}_{\Omega}(u)_r\bigg\rangle =\int_0^t \big\langle u_r\big| x_r\big\rangle d\mu(r). \]
\end{lmma}
The proof follows from Lemma \ref{condition_ito_integral}. \hfill $\Box$
\\
Let $b=(b_t)_{t\in \mathbb{R}_+}$ and $c=(c_t)_{t\in \mathbb{R}_+}$ be centred additive decomposable sections.
Define a measure on $\nu$ on $\mathbb{R}_+$ by $\nu([s,t))=\langle b(s,t)|c(s,t)\rangle$ for every $0\leq s<t.$ Let $x=(x_t)_{t\in \mathbb{R}_+}$ and $y=(y_t)_{t\in \mathbb{R}_+}$ be two continuous adapted processes. Then we have,
\begin{equation}\label{ito_identity_1}
\left\langle \int_{s}^{t} x_r db_r\Bigg| \int_{s}^{t} y_r dc_r\right\rangle=\int_{s}^{t} \langle x_r|y_r \rangle d\nu(r).
 \end{equation}
 We require this result in the following theorem and leave it to the reader for verification.
\begin{thm} \label{bijective_correspondence}
\begin{enumerate}
\item[(1)]  The map  $\text{Exp}_{\Omega}: \mathcal{A}\ni b=(b_t)_{t\in \mathbb{R}_+}\mapsto \text{Exp}(b)=(\text{Exp}(b)_t)_{t\in \mathbb{R}_+}\in \mathcal{C}$ is a bijection.
\item[(2)] Let $b=(b_t)_{t\in \mathbb{R}_+}$ and $c=(c_t)_{t\in \mathbb{R}_+}$ be two centred additive decomposable sections. Then we have
\[ \big\langle \text{Exp}(b)_t\big| \text{Exp}(c)_t\big\rangle =e^{ \langle b_t|c_t \rangle} \text{ for every }t\geq 0.\]
\end{enumerate}
\end{thm}
We remark here that the above Theorem \ref{bijective_correspondence} remains true when we replace the additive decomposible section $\{b_t\}_{t\in \mathbb{R}_+}$ and the corresponding coherent section $\{\text{Exp}(b)_t\}_{t\in \mathbb{R}_+}$ by the additive decomposible family $\{b_t\}_{0<t\leq T}$ and the  coherent family $\{\text{Exp}(b)_t\}_{0<t\leq T}$ for any $T>0$.

Let us recall that the definition of $e$-Logarithm $L^e$. For $t>0$ and $x,y\in D(t)$, we say that $x\sim y$ if there exists a non-zero complex number $\lambda$ such that $x=\lambda y$. Then $\sim$ defines an equivalence relation on $D(t)$. Denote by $\dot{x}$ the equivalence class of $x$ and by $\Delta(t)$ the equivalence classes of $D(t)$. Let $\Delta^{(2)}=\{(t;x,y): x,y\in D(t) \text{ for some} \; t>0\}.$ We say that a function $f:\Delta^{(2)}\to \mathbb{C}$ is continuous if for any given coherent sections $\{x_t\}_{t\in\mathbb{R}_+}$ and $\{y_t\}_{t\in\mathbb{R}_+}$, the map $(0,\infty)\ni t\mapsto f(t;x_t,y_t)\in\mathbb{C}$ is continuous. We say that $f:\Delta^{(2)}\to \mathbb{C}$ vanishing at zero if the limit $\displaystyle \lim_{t\to 0^+} f(t;x_t,y_t)=0$.  Let $e=\{e_t\}_{t\in\mathbb{R}_+}$ be a left coherent decomposible section. By Theorem \cite[Theorem 6.4.2]{Arv}, there exists a unique continuous function $L^e:\Delta^{(2)}\to \mathbb{C}$ vanishing at zero such that 
\[e^{L^e(t;\dot{ x},\dot{ y})}=\frac{ \langle x|y \rangle}{ \langle x|e_t \rangle \langle e_t|y \rangle}.\]
The function $L^e$ is called the $e$-Logarithm. As a consequence of the above theorem we have the following corollory.
\begin{crlre}
 Let $e=\{e_t\}_{t\in\mathbb{R}_+}$ be a centred coherent decomposible section. Then the e-Logarithm is positive definite. More precisely for every $t>0$, the map $D(t)\times D(t)\ni (x,y)\mapsto L^e(t;\dot{x},\dot{y})\in\mathbb{C}$ is positive definite.
\end{crlre}
\begin{prf}
 For $x,y\in D(t)$, let $\{x_s:0<s\leq t\}$ and $\{y_s:0<s\leq t\}$ be the left coherent decomposible families such the $x_t=x$ and $y_t=y$. Then by Theorem  \ref{bijective_correspondence}, there exist additive decomposible families $\{b_s:0<s\leq t\}$ and $\{b_s:0<s\leq t\}$ such that 
 \[\langle x_s| y_s\rangle =e^{ \langle b_s|c_s \rangle} \text{ for every }0<s\leq t.\]
 Set $b=b_t$ and $c=c_t$.
 Recall that $L^e$ is homogeneous i.e. For $t>0$, we have $L(t;\lambda x,\mu y)=L(t; x, y)$ where $x,y\in D(t)$ and $\lambda,\mu\neq 0$.
 
 For $t>0$, let $x,y\in D(t)$. Since $L^e$ is homogeneous, we can assume that $\langle x|e_t\rangle=1$ and $\langle y|e_t\rangle=1$. With the foregoing notation we have $e^{L^e(t;\dot{x},\dot{y})}=\langle x|y\rangle =\langle x_t|y_t\rangle =e^{\langle b_t|c_t\rangle}=e^{\langle b|c\rangle}.$ Hence for $x,y\in D(t)$ there exist unique $b,c\in E(t)$ such that $L^e(t;\dot{x},\dot{y})=\langle b|c\rangle$. This implies that for $t>0$, the map $D(t)\times D(t)\ni (x,y)\mapsto L^e(t;\dot{x},\dot{y})\in\mathbb{C}$ is positive definite.   \hfill $\Box$
\end{prf}
\section{Decomposable vectors of one-parameter CAR flows}
In this section we describe left coherent sections for one-parmeter CAR flows. We achieve this by using the bijective correspondence between the set of all additive decomposable sections and the set of all left coherent sections obtained in the previous section.
Let $H$ be a Hilbert space and let $H^{\otimes n}$ be the n-fold tensor product of $H$ for $n\in \mathbb{N}$. For $\sigma\in S_n$, define a unitary $U_{\sigma}$ on $H^{\otimes n}$ by
\[U_{\sigma}(\xi_1\otimes\xi_2\otimes \cdots\xi_n)=\xi_{\sigma(1)}\otimes\xi_{\sigma(2)}\otimes \cdots\xi_{\sigma(n)}\; \text{ for every }\xi_1,\xi_2,...,\xi_n\in H.\]
Let $H^{\textcircled{\small{a}}^n}$ be the subspace of $H^{\otimes n}$ given by
\[ H^{\textcircled{\small{a}}^n}=\{u\in H^{\otimes n}:U_{\sigma}(u)=\varepsilon(\sigma)u \text{ for all }\sigma\in S_n \}.\] 
Here $\varepsilon(\sigma)$ is 1 if $\sigma$ is even and -1 if $\sigma$ is odd.
We define $\xi_1\wedge\xi_2\wedge \cdots\wedge\xi_n\in H^{\textcircled{\small{a}}^n}$ as
\[\xi_1\wedge\xi_2\wedge \cdots\wedge\xi_n=\frac{1}{\sqrt{n!}}\sum_{\sigma\in S_n}\text{sgn}(\sigma)\xi_{\sigma(1)}\otimes\xi_{\sigma(2)}\otimes \cdots\otimes\xi_{\sigma(n)}\]
and the inner product on $H^{\textcircled{\small{a}}^n}$ as 
\[\langle \xi_1\wedge\xi_2\wedge \cdots\wedge\xi_n|\eta_1\wedge\eta_2\wedge \cdots\wedge\eta_n  \rangle=\text{det}((\xi_i|\eta_j)).\]
Let $\Gamma_a(H)$ be the antisymmetric Fock space given by 
\begin{equation*}
 \Gamma_a(H)=\displaystyle\bigoplus_{n=0}^{\infty} H^{\textcircled{\small{a}}^n}=\mathbb{C}\Omega\oplus \bigoplus_{n=1}^{\infty} H^{\textcircled{\small{a}}^n}.
\end{equation*}
Here $\Omega$ is the vacuum vector.
Let $H_1$ and $ H_2 $ be Hilbert spaces. Then the map $\Gamma_a(H_1)\otimes \Gamma_a(H_2) \to \Gamma_a(H_1\oplus H_2)$ given by
\[(\xi_1\wedge \xi_2\wedge\cdots \wedge \xi_n)\otimes (\eta_1\wedge \eta_2\wedge\cdots \wedge \eta_m)\mapsto \xi_1\wedge \xi_2\wedge\cdots \wedge \xi_n\wedge\eta_1\wedge \eta_2\wedge\cdots \wedge \eta_m  \]
for $\xi_i\in H_1$ and $\eta_j\in H_2$ with $1\leq i\leq n$ and $1\leq j\leq m$, $m,n\in \mathbb{N},$ extends to a unitary operator. We freely use this identification in the rest of the paper.
For $\xi\in H$, define a bounded operator $a^*(\xi)$ on $\Gamma_a(H)$ by
\begin{equation}
 a^*(\xi)\eta:=\begin{cases}
 \xi  & \mbox{ if
} \eta=\Omega,\cr
   
    \xi\wedge \eta &  \mbox{ if } \eta\perp \Omega.
         \end{cases}
\end{equation}
and denote the adjoint of $a^*(\xi)$ by $a(\xi)$. The operators $a^*(\xi)$ and $a(\xi)$ are called the creation and the annihilation operator associated to a vector $\xi$. For an isometric representation $V$ of $P$ on $H$, there exists a unique $E_0$-semigroup $\beta=\{\beta_x\}_{x\in P}$ on $B(\Gamma_a(H))$ satisfying
\[\beta_x(a(\xi))=a(V_x\xi)\]
for every $\xi\in H$, called the CAR flow associated to an isometric representation $V$ of $P$; see \cite{R19}. Let $x\in P$ be given.
Define a unitary $U_x: \Gamma_a(\text{Ker}(V_x^*))\otimes \Gamma_a(H) \to \Gamma_a(H)$ by
\[(\xi_1\wedge \xi_2\wedge\cdots \wedge \xi_n)\otimes (\eta_1\wedge \eta_2\wedge\cdots \wedge \eta_m)\mapsto V_x\eta_1 \wedge V_x\eta_2\wedge\cdots \wedge V_x\eta_m\wedge\xi_1\wedge \xi_2\wedge\cdots \wedge \xi_n  \]
for $\xi_i\in \text{Ker}(V_x^*)$ with $1\leq i\leq n$ and $\eta_j\in H$ with $1\leq j\leq m$, $m,n\in \mathbb{N}$.  With little abuse of notation we write for $\xi\in \Gamma_a(\text{Ker}(V_x^*))$ and $\eta \in  \Gamma_a(H)$, $U_x(\xi\otimes \eta)= \Gamma_a(V_x)\eta\wedge \xi.$
For $x\in P$ and $\xi\in \Gamma_a( \text{Ker}(V_x^*))$, define a bounded operator $T^x_{\xi}:\Gamma_a(H)\to \Gamma_a(H)$ by $T^x_{\xi}\eta=\Gamma_a(V_x)\eta\wedge \xi$. The product system is given by $E(x)=\{T^x_{\xi}: \xi\in \Gamma_a( \text{Ker}(V_x^*))\}.$ For notational convenience, in many places we simply write $\xi$ for $T^x_{\xi}$ in our calculation. 

Let $A$ be a $P$-module and $b\in \Omega$ be given. Define a function $\psi_b^A: \mathbb{R}^d\to \mathbb{R}$ by $\psi_b^A(x)=\text{sup}\{t\in \mathbb{R}: x-tb\in A \}$ for $x\in \mathbb{R}^d$. We simply write $\psi_b$ when $A$ is clear from the context.
 For $k\in\mathbb{N},$ denote the set $\{(r_1,r_2,...,r_k)\in \mathbb{R}^k: r_i=r_j \text{ for some }i\neq j\text{ with } 1\leq i,j\leq d\}$ by $N$ which is a null-set of $\mathbb{R}^k$, and define $\varepsilon^{(k)}:\mathbb{R}^k\to \{-1,0,1\}$ by
\begin{equation}
 \varepsilon^{(k)}(r):=\begin{cases}
 0  & \mbox{ if
} r=(r_1,r_2,...,r_k)\in N,\cr
   
    \text{sgn}(\sigma) &  \mbox{ if } r\notin N \text{ and } \sigma\in S_k \text{ such that } r_{\sigma(1)}>r_{\sigma(2)}>\cdots >r_{\sigma(k)}.
         \end{cases}
\end{equation}
Define a map $ \varepsilon_b^{(k)}:A^k\to \{-1,0,1\}$ by $ \varepsilon_b^{(k)}(x_1,x_2,...,x_k)=\varepsilon^{(k)}(\psi_b(x_1),\psi_b(x_2),...,\psi_b(x_k)) $ for $(x_1,x_2,...,x_k)\in A^k$, and for $\xi\in L^2(A,K)$ define $e^{\varepsilon_b}(\xi)\in \Gamma_a(L^2(A,K))$ by
\[e^{\varepsilon_b}(\xi)=\displaystyle\sum_{k=0}^{\infty} \frac{\varepsilon_b^{(k)}\xi^{\otimes k}}{\sqrt{k!}}\]
Let $A$ be a $P$-module i.e. $A$ is a non-empty closed subset of $\mathbb{R}^d$ such that $A+P\subseteq A$. Let $K$ be a Hilbert space of dimensional $k$ with $k\in \mathbb{N}$. Denote the space of all $K$-valued square integrable functions on $A$ by $L^2(A,K)$. For $x\in P$, define an operator $V_x^{(A,K)}$ on $L^2(A,K)$ by
\[(V_x^{(A,K)}\xi)(y)=\begin{cases}
 f(y-x)  & \mbox{ if
} y-x\in A,\cr
   
    0 &  \mbox{ if }  y-x\notin A.
         \end{cases}\]
Then $V^{(A,K)}$ defines an isometric representation of $P$, called the isometric representation of $P$ associated to $A$ of multiplicity $k$. 
\begin{lmma} \label{Relation_1}
 Let $b,c\in \Omega$ be given. For any $t\geq 0$ and $\eta\in L^2(A,K)$, we have
 \begin{enumerate}
  \item [(1)] $\Gamma_a(V^{(A,K)}_{tb})(\varepsilon_c^{(k)} \eta^{\otimes k})= \varepsilon_c^{(k)} (V^{(A,K)}_{tb}\eta)^{\otimes k}.$
  \item[(2)] $\Gamma_a(V^{(A,K)}_{tb})e^{\varepsilon_c}(\eta)=e^{\varepsilon_c}(V^{(A,K)}_{tb}\eta).$
 \end{enumerate}
\end{lmma}
\begin{prf}
 Let $x_1,x_2,...,x_k\in A$, we have
 \begin{align*}
  (\Gamma_a(V^{(A,K)}_{tb})&\varepsilon_c^{(k)} \eta^{\otimes k})(x_1,x_2,...,x_k) \\
  &= (\varepsilon_c^{(k)} \eta^{\otimes k})(x_1-tb,x_2-tb,...,x_k-tb)1_{A^k}(x_1-tb,x_2-tb,...,x_k-tb)\\
  &= \varepsilon^{(k)}(\psi^A_c( x_1-tb), \psi^A_c(x_2-tb),...,\psi^A_c( x_k-tb)) \\
  &\quad\quad\quad(\eta^{\otimes k})(x_1-tb,x_2-tb,...,x_k-tb)1_{A^k}(x_1-tb,x_2-tb,...,x_k-tb)\\
  &= \varepsilon^{(k)}(\psi^{A+tb}_c( x_1), \psi^{A+tb}_c(x_2),...,\psi^{A+tb}_c( x_k)) (V^{(A,K)}_{tb}\eta)^{\otimes k}(x_1,x_2,...,x_k)\\
  &= \varepsilon^{(k)}(\psi^{A}_c( x_1), \psi^{A}_c(x_2),...,\psi^{A}_c( x_k)) (V^{(A,K)}_{tb}\eta)^{\otimes k}(x_1,x_2,...,x_k)\\
  &=(\varepsilon_c^{(k)} (V^{(A,K)}_{tb}\eta)^{\otimes k})(x_1,x_2,...,x_k).
 \end{align*}
The above equality holds for almost every $(x_1,x_2,...,x_k)\in A^k$. This implies part(1). Clearly part(2) follows from part(1).  
\hfill $\Box$
\end{prf}

Fix $b\in\Omega$. Denote the CAR flow associated to the isometric representation $\{V^{(A,K)}_{tb}\}_{t\geq 0}$ by $\{\beta_t\}_{t\geq 0}$.  
We leave it to the reader to verify that $T^{tb}_{e^{\varepsilon_b}(E^{\perp}_{tb}\xi)}$ is a decomposable vector of $\beta$ for any $t>0$ and $\xi\in L^2(A,K)$. In fact, we have the following proposition.
\begin{ppsn}\label{one_parameter_decomposible}
The set of all decomposable vectors of $\{\beta_t\}_{t\geq 0}$ is given by $\{\lambda T_{e^{\varepsilon_b}(E^{\perp}_{tb}\xi)}: \lambda \in \mathbb{C},\; t\geq 0 \text{ and } \xi\in L^2(A,K)\}$.
\end{ppsn}
\begin{prf}
Let $b=(b_t)_{t\geq 0}$ be an additive decomposable section for $\beta$. Then for $t\geq 0$, $b_t=E^{\perp}_{ta}\xi$ for some $\xi\in L^2(A,K)$. By the proof of Proposition \ref{construction_exponential}, the corresponding left coherent section $u=(u_t)_{t\geq 0}$ is given by \[u_t=\displaystyle \sum_{k=0}^{\infty} x_t^{(k)} \text{ where } x_t^{(0)}=\Omega_t,\;x_t^{(1)}=b_t\;\text{ and } x_t^{(k)}=\int_{0}^{t} x_r^{(k-1)} db_r \;\text{ for }t\geq 0.   \]
 First let us compute $x_t^{(2)}$ for $t\geq 0$. For each $n\in \mathbb{N}$, let $r_{i,n}=\frac{ita}{n}$ with $0\leq i\leq n$.
 \begin{align*}
  x_t^{(2)}&=\int_{0}^{t} E_{ra}^{\perp}\xi dE_{ra}^{\perp}\xi \\
           &= \lim_{n\to \infty} \sum_{i=0}^{n-1} T^{ r_{i,n}}_{E^{\perp}_{r_{i,n}}\xi}\;\circ \; T^{r_{i+1,n}-r_{i,n}}_{V^*_{r_{i,n}}E^{\perp}_{r_{i+1,n}}\xi}\;\circ\; T_{\Omega_{ta-r_{i+1,n}}}^{ta-r_{i+1,n}} \;(\text{ by Proposition }\ref{condition_ito_integral}) \\
           &=\displaystyle \lim_{n\to \infty}  T^{ta}_{\sum_{i=0}^{n-1} E_{r_{i,n}}E^{\perp}_{r_{i+1,n}}\xi\wedge E^{\perp}_{r_{i,n}}\xi }\\
           &=\displaystyle \lim_{n\to \infty}  \sum_{i=0}^{n-1} E_{r_{i,n}}E^{\perp}_{r_{i+1,n}}\xi\wedge E^{\perp}_{r_{i,n}}\xi .
           \end{align*}
In the view of the Proposition \ref{condition_ito_integral} it is enough to check the pointwise convergence of $\sum_{i=0}^{n-1} E_{r_{i,n}}E^{\perp}_{r_{i+1,n}}\xi\wedge E^{\perp}_{r_{i,n}}\xi$ to $\frac{\varepsilon^{(2)}_a \; (E_{ta}^{\perp}\xi)^{\otimes 2}}{\sqrt{2!}}$ almost everywhere.           
For $x\in A$, there exists a unique $\widetilde{x}\in \partial A$ such that $x=\widetilde{x}+\psi_a(x)a$. For almost every $(x,y)\in A\times A$ and for large $n$, we have
\begin{align*}
  &\bigg(  \sum_{i=0}^{n-1} E_{r_{i,n}}E^{\perp}_{r_{i+1,n}}\xi\wedge E^{\perp}_{r_{i,n}}\xi\bigg)(x,y)
 = \frac{1}{\sqrt{2!}}\sum_{i=0}^{n-1}\big(\chi_{(A+\frac{ita}{n})\cap (A\setminus A+\frac{(i+1)ta}{n})}(x)  \chi_{(A\setminus A+\frac{ita}{n})}(y) \\
 &\quad\quad\quad\quad-   \chi_{(A\setminus A+\frac{ita}{n})}(x) \chi_{(A+\frac{ita}{n})\cap (A\setminus A+\frac{(i+1)ta}{n})}(y)\big) \xi(x)\otimes \xi(y)     
 \end{align*}
When $\psi_a(x)<\psi_a(y)$, $x\in A\setminus A+\frac{ita}{n}$ and $y\in (A+\frac{ita}{n})\cap (A\setminus A+\frac{(i+1)ta}{n})$ for some $0\leq i\leq n-1$. Hence $\bigg(  \sum_{i=0}^{n-1} E_{r_{i,n}a}E^{\perp}_{r_{i+1,n}a}\xi\wedge E^{\perp}_{r_{i,n}a}\xi\bigg)(x,y)=\frac{-(E_{ta}^{\perp}\xi)(x)\otimes (E_{ta}^{\perp}\xi)(y)}{\sqrt{2!}}$. Similarly if $\psi_a(x)>\psi_a(y)$, $\bigg(  \sum_{i=0}^{n-1} E_{r_{i,n}a}E^{\perp}_{r_{i+1,n}a}\xi\wedge E^{\perp}_{r_{i,n}a}\xi\bigg)(x,y)=\frac{(E_{ta}^{\perp}\xi)(x)\otimes (E_{ta}^{\perp}\xi)(y)}{\sqrt{2!}}$. We conclude that $x_t^{(2)}=\int_{0}^{t} E_{ra}^{\perp}\xi dE_{ra}^{\perp}\xi=\frac{\varepsilon^{(2)}_a \; (E_{ta}^{\perp}\xi)^{\otimes 2}}{\sqrt{2!}}$. 
 
Before proving $x_t^{(k)}=\frac{\varepsilon^{(k)}_a \; (E_{ta}^{\perp}\xi)^{\otimes k}}{\sqrt{k!}}$ for any $k\in\mathbb{N}$. Let us fix few notation. For $\xi_1,\xi_2,...,\xi_n,\eta\in H$, set $\eta^{(k)}\odot (\xi_1\otimes \xi_2\otimes\cdots\xi_{n})=\xi_1\otimes \xi_2\otimes\cdots \xi_{k-1}\otimes \eta\otimes \xi_{k}\otimes\cdots\xi_{n}$ for $1\leq k\leq n$.  With this notation we can see that
\[\xi_1\wedge (\xi_2\wedge \xi_3\wedge\cdots\wedge\xi_n)=\displaystyle\frac{1}{\sqrt{k}}\sum_{j=1}^n (-1)^{j-1} \xi_1^{(j)}\odot (\xi_2\wedge \xi_3\wedge\cdots\wedge\xi_n).\]
Assume that $x_t^{(k-1)}=\frac{\varepsilon^{(k-1)}_a \; (E_{ta}^{\perp}\xi)^{\otimes (k-1)}}{\sqrt{k-1!}}$ for some $k\in \mathbb{N}$.
Now consider the following expression. For almost every $(x_1,x_2,...,x_k)\in A^k$ and for large $n$, we have
\begin{align*}
\big( \sum_{i=0}^{n-1}E_{r_{i,n}}&E^{\perp}_{r_{i+1,n}}\xi\wedge  \frac{\varepsilon^{(k-1)}_a \; (E_{r_{i,n}a}^{\perp}\xi)^{\otimes k-1}}{\sqrt{(k-1)!}} \big)(x_1,x_2,...,x_k)\\
&=\sum_{i=0}^{n-1}\frac{1}{\sqrt{k}}\sum_{j=1}^k (-1)^{j-1} E_{r_{i,n}}E^{\perp}_{r_{i+1,n}}\xi\odot  \frac{\varepsilon^{(k-1)}_a \; (E_{r_{i,n}a}^{\perp}\xi)^{\otimes k-1}}{\sqrt{(k-1)!}} \big)(x_1,x_2,...,x_k)\\ 
&=\sum_{i=0}^{n-1}\frac{1}{\sqrt{k!}}\sum_{j=1}^k  (-1)^{j-1}
\chi_{(A+\frac{ita}{n})\cap (A\setminus A+\frac{(i+1)ta}{n})}(x_j) \varepsilon^{(k-1)}_a(x_1,x_2,...\widehat{x_j},...,x_k) 
\\
&\quad\quad\prod_{l=1,l\neq j}^{k}  \chi_{(A+\frac{ita}{n})}(x_l) (E_{ta}^{\perp}\xi)^{\otimes k} (x_1,x_2,...,x_k).
\end{align*}
There exist unique $i$ and $j$ such that $x_j\in(A+\frac{ita}{n})\cap (A\setminus A+\frac{(i+1)ta}{n})$ and $x_1,x_2,...\widehat{x_j},...,x_k\in A+\frac{ita}{n}$.
Hence $\big( \sum_{i=0}^{n-1}E_{r_{i,n}}E^{\perp}_{r_{i+1,n}}\xi\wedge  \frac{\varepsilon^{(k-1)}_a \; (E_{r_{i,n}a}^{\perp}\xi)^{\otimes k-1}}{\sqrt{(k-1)!}} \big)=\frac{\varepsilon^{(k)}_a \; (E_{ta}^{\perp}\xi)^{\otimes k}}{\sqrt{k!}}$ almost everywhere.
By definition we have
\begin{equation*}
 x_t^{(k)}=\int_{0}^{t}   \frac{\varepsilon^{(k-1)}_a \; (E_{ra}^{\perp}\xi)^{\otimes k-1}}{\sqrt{(k-1)!}}   dE_{ra}^{\perp}\xi 
 =\displaystyle \lim_{n\to \infty}  \sum_{i=0}^{n-1} E_{r_{i,n}}E^{\perp}_{r_{i+1,n}}\xi\wedge  \frac{\varepsilon^{(k-1)}_a \; (E_{r_{i,n}a}^{\perp}\xi)^{\otimes k-1}}{\sqrt{(k-1)!}}.
\end{equation*}
Hence we conclude that for each $k\in \mathbb{N}$, $x_t^{(k)}=\frac{1}{\sqrt{k!}}\varepsilon_a^{(k)}(E_{ta}^{\perp}\xi)^{\otimes k}$ and
$u_t=e^{\varepsilon_a}(E_{ta}^{\perp}\xi)$.
The proof is complete.
\hfill $\Box$
\end{prf}

\section{Characterization for Decomposablility of CAR flows} 
Let us recall the definition of decomposable product system from \cite{sundar2019arvesons}. Let $\alpha=\{\alpha_x\}_{x\in P}$ be an $E_0$-semigroup over $P$ on $B(H)$. For $x\in P$, let $E(x)=\{T\in B(H): \alpha_x(A)T=TA \text{ for all }A\in B(H)\}.$ Then $E=\{E(x):x\in P\}$ has the structure of the product system.

For $x\in P$, a non-zero element $u\in E(x)$ is said to be a decomposable vector if for $y\in P$ with $y\leq x$, then there exist $v\in E(y)$ and $w\in E(x-y)$ such that $u=vw$. Denote the set of all decomposable vectors in $E(x)$ by $D(x)$. We say that the product system $E=\{E(x):x\in P\}$ is decomposable if the following conditions are satisfied.
\begin{enumerate}
 \item For each $x,y\in P$, $D(x)D(y)\subseteq D(x+y)$.
 \item For each $x\in P$, $D(x)$ is total in $E(x)$.
\end{enumerate}
We call an $E_0$-semigroup is decomposable if its associated product system is decomposable.
Let $A$ be a $P$-module and let $K$ be a Hilbert space of dimension $k$ with $k\in \mathbb{N}$. 
Let us denote the CAR flow associated to the isometric representation $V^{(A,K)}$ by $\beta$. The goal of this section is to exhibit the necessary and sufficient condition for the CAR flow $\beta$ to be decomposable. With the foregoing notation we have the following proposition.
\begin{ppsn} \label{epsilon_equality}
Let $E$ be the product system associated to $\beta$. Assume that $E=\{E(x): x\in P\}$ is a decomposable product system. Then for any given $b,c\in \Omega$, we have
\[\varepsilon^{(2)}_b(E_{sb}^{\perp}E_{tc}^{\perp}\xi)^{\otimes 2}=\varepsilon^{(2)}_c(E_{sb}^{\perp}E_{tc}^{\perp}\xi)^{\otimes 2} \;a.e \text{ for all } \xi\in L^2(A,K).\]
\end{ppsn}
\begin{prf}
 Without loss of generality we assume that $b$ and $c$ are linearly independent. Choose a linearly independent collection $v_1,v_2,...,v_d\in\Omega$  with $v_1=b$ and $v_2=c$. Let $Q=\{r_1v_1+r_2v_2+\cdots+r_dv_d:\text{ for all } r_i\geq 0 \text{ with }1\leq i\leq d\}$.
 Since $E$ is decomposable over $P$, it is also decomposable over $Q$. For the product system $\{E(x)\}_{x\in Q}$ over $Q$, let $D(x)$ be the set of all decomposable vectors in $E(x)$.
 By Proposition \ref{one_parameter_decomposible}, we have
 \begin{align*}
  D(sb)&=\{\lambda T_{e^{\varepsilon_b}(E^{\perp}_{sb}\xi)}: \lambda \in \mathbb{C},\; s\geq 0 \text{ and } \xi\in L^2(A,K)\}, \text{ and}\\
  D(tc)&=\{\mu T_{e^{\varepsilon_c}(E^{\perp}_{tc}\eta)}: \mu \in \mathbb{C},\; t\geq 0 \text{ and } \eta\in L^2(A,K)\}.
 \end{align*}
As $E$ is decomposable over $Q$, $D(sb)D(tc)=D(sb+tc)=D(tc)D(sb)$ i.e. for any $\xi, \eta \in L^2(A,K)$, there exist $\xi', \eta' \in L^2(A,K)$ such that 
\[T^{sb}_{e^{\varepsilon_b}(E^{\perp}_{sb}\xi)}T^{tc}_{e^{\varepsilon_c}(E^{\perp}_{tc}\eta)}=T^{tc}_{e^{\varepsilon_c}(E^{\perp}_{tc}\eta')}T^{sb}_{e^{\varepsilon_b}(E^{\perp}_{sb}\xi')}.\]
By applying $\zeta$ on both sides, we have 
\[\Gamma_a(V_{sb+tc})\zeta\wedge\Gamma_a(V_{sb})e^{\varepsilon_c}(E^{\perp}_{tc}\eta)\wedge e^{\varepsilon_b}(E^{\perp}_{sb}\xi)=
 \Gamma_a(V_{sb+tc})\zeta\wedge \Gamma_a(V_{tc})e^{\varepsilon_b}(E^{\perp}_{sb}\xi')\wedge e^{\varepsilon_c}(E^{\perp}_{tc}\eta').\] 
 The above equation together with Lemma \ref{Relation_1}, 
 \[e^{\varepsilon_c}(V_{sb}E^{\perp}_{tc}\eta)\wedge e^{\varepsilon_b}(E^{\perp}_{sb}\xi)=
  e^{\varepsilon_b}(V_{tc}E^{\perp}_{sb}\xi')\wedge e^{\varepsilon_c}(E^{\perp}_{tc}\eta').\]
Applying $\Gamma_a(E^{\perp}_{sb}E^{\perp}_{tc})$ on both sides, $e^{\varepsilon_b}(E^{\perp}_{sb}E^{\perp}_{tc}\xi)= e^{\varepsilon_c}(E^{\perp}_{sb}E^{\perp}_{tc}\eta').$ Equating $1$-particle space and $2$-particle space on both sides, we have $E^{\perp}_{sb}E^{\perp}_{tc}\xi=E^{\perp}_{sb}E^{\perp}_{tc}\eta'$ and hence 
\[\varepsilon^{(2)}_b(E_{sb}^{\perp}E_{tc}^{\perp}\xi)^{\otimes 2}=\varepsilon^{(2)}_c(E_{sb}^{\perp}E_{tc}^{\perp}\xi)^{\otimes 2} \;a.e \text{ for all } \xi\in L^2(A,K). \]
\hfill $\Box$
\end{prf}
\begin{ppsn} \cite[Proposition 2.3(3)]{AS19} \label{interior_A}
 Let $A$ be a $P$-module and let $a\in \Omega$ be given. Then the map $\partial A\times (0,\infty)\ni (x,t)\mapsto x+ta\in \text{Int}(A)$ is a homeomorphism. 
\end{ppsn}
Let $v_1,v_2,...,v_d\in \Omega$ be a linearly independent set in $\mathbb{R}^d$. Define a function $\varphi: \mathbb{R}^{d-1}\to \mathbb{R}$ by $\varphi(r_1,r_2,...,r_{d-1}):=r_d-\psi_{v_d}(r_1v_1+r_2v_2+\cdots+r_dv_d)$. The map is well-defined. This follows from the observation that for $s\in \mathbb{R}$, $\psi_{v_d}(r_1v_1+r_1v_1+\cdots (r_{d}+s)v_d)=s+\psi_{v_d}(r_1v_1+r_1v_1+\cdots r_{d}v_d)$. 
 With the foregoing notation we have the following lemma. 
\begin{lmma}\label{boundary_A}Let $a\in \Omega$ be given. Then we have the following.
\begin{enumerate}
 \item [(1)] The map $\partial A\times \mathbb{R}\ni (x,t)\mapsto x+ta \in \mathbb{R}^d$ defines a homeomorphism.
 \item[(2)] The $P$-module $A$ is given by $A=\{x\in \mathbb{R}^d: \psi_a(x)\geq 0\}$.
 \item[(3)] The boundary of $A$ is given by $\partial A=\{x\in \mathbb{R}^d: \psi_{v_d}(x)=0\}$.
 Moreover, the map $B:\mathbb{R}^{d-1}\to \partial A$ given by 
 \[B(r_1,r_2,...,r_{d-1}):=r_1v_1+r_2v_2+\cdots+r_{d-1}v_{d-1}+\varphi(r_1,r_2,...,r_{d-1})v_d\]
 for $(r_1,r_2,...,r_{d-1})\in \mathbb{R}^{d-1}$, defines a 
 homeomorphism.
 \end{enumerate}
 \end{lmma}
 \begin{prf}
The proof of part(1) follows similar to \cite[Proposition 2.3(3)]{AS19}. Part(2) is clear from the definition. Note that $\partial A=\{x\in \mathbb{R}^d: \psi_{v_d}(x)= 0\}$.
Define $\varphi:\mathbb{R}^{d-1}\to \mathbb{R}$ by
\[\varphi(r_1,r_2,...,r_{d-1})=r_d-\psi_{v_d}(r_1v_1+r_2v_2+\cdots+r_{d-1}v_{d-1}+ r_dv_d)\]
for every $(r_1,r_2,...,r_{d-1},r_d)\in\mathbb{R}^d$. Now the part(3) is clear from the fact that $\partial A=\{x\in \mathbb{R}^d: \psi_{v_d}(x)= 0\}$ and the remaining we leave it to the reader for verification.
\hfill $\Box$
 \end{prf}
\begin{thm}\label{decomposible_implies_hyperplane}
 Let $E$ be the product system over $P$ corresponding to $\beta$. Suppose $E$ is decomposable, then there exists an element $\lambda\in P^*$ such that $A=\{x\in \mathbb{R}^d: \langle x|\lambda\rangle \geq 0\}$.
\end{thm}
\begin{prf}
Let $v_1,v_2,...,v_d\in \Omega$ be a linearly independent set in $\mathbb{R}^d$.
 Fix $i$ with $1\leq i\leq d-1$ and $t_0>0$. By Proposition \ref{interior_A} and Lemma \ref{boundary_A}(3), there exist functions $f_0:\mathbb{R}^{d-1}\to \mathbb{R}^{d-1}$ and $g_0:\mathbb{R}^{d-1}\to (0,\infty)$ such that for every $r=(r_1,r_2,...,r_{d-1})\in \mathbb{R}^{d-1}$, we have
 \[B(r)+t_0v_i=B(f_0(r))+g_0(r)v_d\]
 i.e $r_1v_1+r_2v_2+\cdots+r_{d-1}v_{d-1}+\varphi(r)v_d+t_0v_i=i_1\circ f_0(r)v_1+i_2\circ f_0(r)v_2+\cdots+i_{d-1}\circ f_0(r)v_{d-1}+\varphi(f_0(r))v_d+g_0(r)v_d$. Here $i_l:\mathbb{R}^{d-1}\to \mathbb{R}$ denote the projection onto the $l^{th}$ coordinate for each $1\leq l\leq d-1$.
 By equating the coefficients of $v_l'$s on both sides, we see that $f_0(r)=r+t_0e_i$ and $g_0(r)=\varphi(r)-\varphi(r+t_0e_i)$. Here $e_i=(0,0,..,1,..0)\in \mathbb{R}^{d-1}$, where $1$ in the $i^{th}$ place and elsewhere $0$. 
 Hence $f_0$ and $g_0$ are continuous. Similarly for $\delta>0$, there exist continuous functions $f_{\delta}:\mathbb{R}^{d-1}\to \mathbb{R}^{d-1}$ and $g_{\delta}:\mathbb{R}^{d-1}\to (0,\infty)$ such that
 \[B(r)+(t_0+\delta)v_i=B(f_{\delta}(r))+g_{\delta}(r)v_d \text{ for all }r\in \mathbb{R}^{d-1}. \]
 We claim that for each $t_0>0$, $g_0(.)$ is constant. Suppose not, for some $t_0>0$, $g_0(.)$ is not constant. Choose $r', r''\in \mathbb{R}^{d-1}$ such that $g_0(r')>g_0(r'')$. Note that the map $[0,\infty)\times \mathbb{R}^{d-1} \ni (\delta, r)\mapsto g_{\delta}(r)\in (0,\infty)$ is continuous.  Choose $\delta_0>0$ such that $g_0(r')>g_{\delta_0}(r'') $. 
 Then there exists a compact neighbourhood $K$ of $\text{int}(A)\times \text{int}(A)$ containing $\big((t_0,r'),(t_0+\delta_0,r'')\big)$ such that 
 \[t<t' \text{ and }g_{t}(r)>g_{t'}(s),\text{ for every }(B(r)+tv_i,B(s)+t'v_i)\in K.\]
 This implies that $\varepsilon^{(2)}_{v_i}\neq \varepsilon^{(2)}_{v_d}$ on $K$. This is a contradiction to Proposition \ref{epsilon_equality}.  Hence for each $s_0>0$, $g_0(.)$ is a constant function. 

Without loss of generality, we assume that $\partial (A)\ni 0$. Then $\varphi(0)=0$. Using this equality and the fact that $g_0(.)$ is a constant function, we see that for every $i$ with $1\leq i\leq d-1$, $\varphi(r+s_0e_i)=\varphi(r)+\varphi(s_0e_i)$ for $r\in \mathbb{R}^{d-1}$ and $s_0>0$. 
Consequently, we deduce that $\varphi(r)=\varphi(r+c)-\varphi(c)$ for every $r\in \mathbb{R}^{d-1}$ and $c=(c_1,c_2,..,c_{d-1})$ with $c_i>0$.
Let $r,r'\in \mathbb{R}^{d-1}$ be given. Choose $c=(c_1,c_2,..,c_{d-1}),c'=(c_1',c_2',..,c'_{d-1})\in \mathbb{R}^{d-1}$ with $c_i,c_i'>0$  and $ c_i+r_i,c_i'+r_i'>0$  for each $1\leq i\leq d-1$.
\begin{align*}
 \varphi(r)+\varphi(r')&=\varphi(r+c)-\varphi(c)+\varphi(r'+c')-\varphi(c')\\
                    &=\varphi(r+c)+\varphi(r'+c')-(\varphi(c)+\varphi(c'))\\
                    &=\varphi(r+r'+c+c')-\varphi(c+c')=\varphi(r+r').
\end{align*}
Since $\varphi$ is continuous, there exists a unique $\mu\in \mathbb{R}^{d-1}$ such that $\varphi(r)=\langle r|\mu\rangle$ for every $r\in \mathbb{R}^{d-1}$.
By Lemma \ref{boundary_A}(2), $r_d-\varphi(r_1,r_2,...,r_{d-1})=r_d-\langle r|\mu\rangle=\psi_{v_d}(r_1v_1+r_2v_2+\cdots+r_dv_d)\geq 0$ for the points of $A$. Therefore $A=\{x\in \mathbb{R}^d: \langle x|\lambda\rangle \geq 0\}$ for some $\lambda\in \mathbb{R}^{d}$. Since $P\subseteq A$, $\lambda\in P^*$. This completes the proof.  
\hfill $\Box$
\end{prf}
\begin{ppsn} \label{decomposible_productsystem}
 Let $v_1,v_2,...,v_d$ be a linearly independent set in $\mathbb{R}^d$ and let $Q:=\{r_1v_1+r_2v_2+\cdots+r_dv_d: \text{ each }r_i\geq 0\}$. Let $E=\{E(x)\}_{x\in Q}$ be a product system over $Q$. For $1\leq i\leq d$, set $E^{(i)}=\{E^{(i)}(r)=E(rv_i): r\geq 0\}$ and denote the set of all decomposable vectors in $E^{(i)}(r)$ by $D^{(i)}(r)$ for $r\geq 0$. Suppose for $1\leq i,j\leq d$ and $r,s\geq 0$, $D^{(i)}(r)D^{(j)}(s)=D^{(j)}(s)D^{(i)}(r)$. Then $E$ is decomposable over $Q$ if and only if each $E^{(i)}$ is decomposable over $\mathbb{R}_+$.
\end{ppsn}
\begin{prf}
 We only prove for $d=2$ and the proof for general $d$ is similar.
 Assume that $E$ is a decomposable product system. Observe that $D^{(i)}(r)=D(rv_i)$ for any $r\geq 0$ and $E^{(i)}$ is a decomposable product system for $i=1,2$. It follows that $D^{(1)}(r)D^{(2)}(s)=D^{(2)}(s)D^{(1)}(r)$.
 Conversely assume that  each $E^{(i)}$ is decomposable over $\mathbb{R}_+$ and $D^{(1)}(r)D^{(2)}(s)=D^{(2)}(s)D^{(1)}(r)$. We claim that $D(rv_1+sv_2)= D^{(1)}(r)D^{(2)}(s)$.
 Clearly $D(rv_1+sv_2)\subseteq D^{(1)}(r)D^{(2)}(s)$. On the other hand, let $u\in D^{(1)}(r)$ and $v\in D^{(2)}(s)$. Let $(r'-r)v_1+(s'-s)v_2\in Q$ i.e. $r'\geq r$ and $s'\geq s$. Write $u=u'u''\in E^{(1)}(r)E^{(2)}(r'-r)$ and $v=v'v''\in E^{(1)}(s)E^{(2)}(s'-s)$. By using the relation $D^{(1)}(r'-r)D^{(2)}(s)=D^{(2)}(s)D^{(1)}(r'-r)$ we have $u''v'=\tilde{v}\tilde{u}$ for some $\tilde{v}\in D^{(2)}(s)$  $\tilde{u}\in D^{(1)}(r'-r)$ and 
 \[uv=u'u''v'v''=u\tilde{v}\tilde{u}v''\in E(rv_1+sv_2)E((r'-r)v_1+(s'-s)v_2).\]
 Hence $uv\in D(rv_1+sv_2)$. Clearly $D(rv_1+sv_2)$ is total in $E(rv_1+sv_2)$. Now
 \begin{align*}
  D(rv_1+sv_2)D(r'v_1+s'v_2)&=D^{(1)}(r)D^{(2)}(s)D^{(1)}(r')D^{(2)}(s')\\
                            &=D^{(1)}(r)D^{(1)}(r')D^{(2)}(s)D^{(2)}(s')\\
                            &=D^{(1)}(r+r')D^{(2)}(s+s').
 \end{align*}
The proof is complete.
\hfill $\Box$
 \end{prf}
\begin{ppsn} \label{hyperplane_implies_decomposible}
 Let $A=\{x\in \mathbb{R}^d: \langle x|\lambda\rangle \geq 0\}$ for some $\lambda\in P^*$ and let $E$ be the product system for the CAR flow over $P$ associated to the $P$-module $A$ of multiplicity $k$. Then $E$ is decomposable over $P$.
\end{ppsn}
\begin{prf} 
 Since $P\subseteq A$, choose a linearly independent set  $\{v_1,v_2,...,v_d\}$ in $A$ such that $v_1,v_2,...,v_{d-1}\in \partial A$ and $P\subseteq Q=\{r_1v_1+r_2v_2+\cdots+r_dv_d: \text{each }r_i\geq 0\}$. Note that $\partial A$ is a $d-1$ vector space over $\mathbb{R}$.
 Observe that $A$ is also a $Q$-module. Let $F=\{F(x)\}_{x\in Q}$ be the product system for the CAR flow over $Q$ associated to $A$ of multiplicity $k$. Note that $E=\{F(x)\}_{x\in P}$ and
 for $1\leq i\leq d-1$,  $D^{(i)}(r)=F^{(i)}(r)=\{\lambda \Gamma_a(V_{rv_i}): \lambda\in \mathbb{C}\}$. By Proposition \ref{one_parameter_decomposible}, $D^{(d)}(r)=\{\lambda T^{rv_d}_{e^{\varepsilon_d}(\xi)}: \lambda\in \mathbb{C}\text{ and } \xi\in \text{Ker}(V_{rv_d}^*)\}$. Clearly $D^{(i)}(r)D^{(j)}(s)=D^{(j)}(s)D^{(i)}(r)$ for each $1\leq i,j\leq d$ and $r,s\geq 0$. By Proposition \ref{decomposible_productsystem}, $F$ is decomposable over $Q$. This implies that $E$ is decomposable over $P$.   
 \hfill $\Box$
\end{prf}
\begin{thm}
 Let $A$ be a $P$-module and let $E$ be the product system over $P$ for the CAR flow associated to $A$ of multiplicity $k$. Then $E$ is decomposable over $P$ if and only if $A=\{x\in \mathbb{R}^d: \langle x|\lambda\rangle \geq 0\}$ for some $\lambda\in P^*$. 
\end{thm}
\begin{prf}
 Proof follows from the Theorem \ref{decomposible_implies_hyperplane} and the Proposition \ref{hyperplane_implies_decomposible}.
 \hfill $\Box$
\end{prf}

\begin{crlre}
 There are uncountable many CCR flows cocycle conjugate to CAR flows over $P$.
\end{crlre}
\begin{prf}
Let $A=\{x\in \mathbb{R}^d: \langle x|\lambda\rangle \geq 0\}$ and $B=\{x\in \mathbb{R}^d: \langle x|\mu\rangle \geq 0\}$ for some $\lambda,\mu \in P^*$. Let $\alpha$ be the CCR flow of $A$ with multiplicity $k$ and let $\beta$ be the CAR flow of $B$ with multiplicity $l$.

We claim that $\alpha$ is cocycle conjugate to $\beta$ if and only if $A=B$ and $k=l$.
Let $E$ and $F$ be the product systems associated to $\alpha$ and $\beta$ respectively. Assume that $\alpha$ is cocycle conjugate to $\beta$. Then by \cite[Theorem 2.9]{MS}, $E$ is isomorphic to $F$ as product systems. Observe that the product systems $E$ and $F$ are embeddable and its corresponding isometric representations of $P$ are unitary equivalent to $V^{(A,K)}$ and $V^{(B,L)}$  respectively; see \cite[Definition 3.3]{R19} and discussion following that. Since $E$ is isomorphic to $F$, $V^{(A,K)}$ is unitary equivalent to $V^{(B,L)}$. Hence $A=B$ and $k=l$.
Conversely let $A=B$ and $k=l$.
Since $F$ is decomposable and has a unit, by \cite[Theorem 4.4]{sundar2019arvesons} there exists a CCR flow $\alpha^V$ given by an isometric representation $V$ such that its product system is isomorphic to $F$. Since the isometric representation constructed out of $F$ is $V^{(A,k)}$, $V$ is unitary equivalent to $V^{(A,k)}$. Hence $E$ is isomorphic to $F$.
Since there are uncountable many $P$-modules of the form $A=\{x\in \mathbb{R}^d: \langle x|\lambda\rangle \geq 0\}$ for $\lambda\in P^*$, the proof is complete.
\hfill $\Box$
\end{prf}
\section*{Acknowledgment}
I would like to thank S. Sundar for suggesting the project and the useful discussions. The author is supported by The Institute of Mathematical Sciences Postdoctoral fellowship.
                                                                                                                                                                            \bibliography{reference}
 \bibliographystyle{amsplain}

\end{document}